\documentclass[11pt]{amsart}

\usepackage{amsmath,amssymb,amsfonts}

\usepackage{caption2}
\usepackage[]{hyperref}

\usepackage[pdftex]{color,graphicx}
\usepackage{multicol}
\usepackage{verbatim}

\numberwithin{equation}{section}

\setcaptionmargin{1cm}

\newcommand{\be}{\begin{equation}}
\newcommand{\ee}{\end{equation}}

\usepackage{color}
\definecolor{rosso}{cmyk}{0,1,1,0.4}
\definecolor{rossos}{cmyk}{0,1,1,0.55}
\definecolor{rossoc}{cmyk}{0,1,1,0.2}

\begin{document}

\begin{center}
{
 \bf 
About zero counting of Riemann Z function
  }  
\vskip 0.2cm {\large

Giovanni  Lodone \footnote{
giolodone3@gmail.com }  
.} 
\\[0.2cm]

{\bf

ABSTRACT    } 
\\[0.2cm]

\vskip 0.2cm {\large 

 An approximate formula for  complex Riemann Xi function, %
 previously developed in \cite {Giovanni Lodone 2021},  is used to   refine Backlund's estimate of %
 the number of zeros till 
a chosen imaginary coordinate } 
\\[0.2cm]

\vskip 0.1cm

\end{center}

{

MSC-Class  :   11M06,  11M99

  }


                            \tableofcontents

\vspace{0.1cm}

\section{%
 Introduction}

  \noindent As well known Rieman $\zeta(s)$ function, 
 (see also  \cite  {Apostol:1976}, \cite  {Davenport1980}, 
\cite[p.~229] {Edwards:1974cz}), is defined as the analytic continuation of the function:

\be \label {1p1}
\zeta(s)=\sum _{n=1}^\infty \frac{1}{n^s} = \prod_{\forall p} \frac{1}{1-\frac{1}{p^s}} \quad \mbox{ ($p$ prime)}  \quad \quad \Re(s)>1
\ee
The strip
  $0<\Re(s)<1 $ is called critical strip, and,  the line $\Re(s)= 1/2$ is  the critical line. The complex argument $s$ is expressed throughout as:
\be \label{1p2}
 s=\frac{ 1}{2}+\epsilon +it
 \ee
so that $\epsilon$ is the distance from critical line ,%
 and, $t$ is the imaginary coordinate.
 
\noindent A key function is $\xi(s)$ 
  \cite[p.~16]{Edwards:1974cz}, that is real on the critical line and has the same zeros of $\zeta(s)$ in the critical strip:
  
\be \label {1p3}
\xi(s)=\Gamma \left( \frac{s}{2}  +1\right)(s-1)\pi^{-s/2}\zeta(s) = 
\left[ \frac{\Gamma \left( \frac{s}{2}  \right) }{ \pi^{s/2}} \right]  \ \times \   \left[  \zeta(s)(s-1) \ \frac{s}{2}  \right]
\ee

\noindent In \cite {Giovanni Lodone 2021}  is developed  an extension  of the Riemann-Siegel function  \cite[p.~136]{Edwards:1974cz} %
that provides a means of expressing  $\xi(\frac{ 1}{2}+\epsilon +it)$ at least   for $ 0 \le \epsilon  < 2$ and $|t|>>1$ , based on   \cite[~ chp. 7] {Edwards:1974cz}.

\noindent For ease of reading,  we present  final result  
of  \cite {Giovanni Lodone 2021}:

\be \label {ZSinhECosh}
 Z(t,\epsilon)=
2 \sum_{n=1}^N \frac{    \cosh \left[ \epsilon \  \ln\left( \sqrt{\frac {t}{2 \pi n^2}  }\right) \right]   }{\sqrt{n}}
\cos \left(   t \  \ln\left(    \sqrt{\frac {t}{2  e \pi n^2} } \right) - \frac{\pi}{8}   \right) + 
\ee

$$  + 2 i  \sum_{n=1}^N \frac{  \sinh\left[ \epsilon \ \ln\left( \sqrt{\frac {t}{2 \pi n^2}  }\right) \right]   }{\sqrt{n}}
\sin \left(   t  \  \ln\left(    \sqrt{\frac {t}{2  e \pi n^2} } \right) - \frac{\pi}{8}  \right )         +R(t,\epsilon)
+ Err(t,\epsilon)  = %
\frac{-\xi( 1/2+\epsilon +it)}{F(t)  e^{i\epsilon\frac{\pi}{4}}}
$$

\noindent where %
$N=  \left  \lfloor   \sqrt{\frac {t}{2   \pi } }  \right \rfloor $,
$F(t) =  \left(\pi/2 \right)^{0.25}t^{\frac{7}{4}} e^{- \frac{\pi}{4}  t}  \ $,  \cite[p.~119]{Edwards:1974cz} for $t$ big, $ |Err(t,\epsilon)| < e^{-0.1 \times t}$ at least for $|t|>100$  ( \cite[p.~144]{Edwards:1974cz},\cite[p.~37] {Giovanni Lodone 2021}), and  $R_M(t,\epsilon)$  is given by eq. 
(\ref{RDue}), 
 specialized for $M=2$,  (i.e.  $R(t,\epsilon)=R_{M \rightarrow \infty}(t,\epsilon)$ )%
 :

\be \label {RDue}
 R_2(t,\epsilon) = (-1)^{N-1} \left( \frac{2 \pi}{t}\right)^{1/4}  \left[     C_0(p) \omega^0+C_1(p,\epsilon) \omega^1 + C_2(p,\epsilon) \omega^2 \right] \ \ ; \ \ M=2
\ee

\noindent where : 

$$  p= \sqrt{\frac {t}{2 \pi}} -N
\ ; \ \ \ \
 \omega = \left( \frac{2 \pi}{t}\right)^{1/2} \ ; \ \ \ \
 C_0(p)=\psi(p)^{(0)}= \cos(2 \pi(p^2-p-1/16) /  \cos(2 \pi p) 
 $$ besides $\psi(p)^{(n)}$ means the $n^{th}$ derivative with respect to $p$.

\be \label {Cuno}
C_1(p,\epsilon) 
 =  -  \frac{ i
  \psi^{(1)}(p)}{4 \pi}    \epsilon 
-  \frac{\psi^{(3)}(p)}{2^5 \pi^2 3} 
\ee

\noindent Notice that, for $\epsilon=0$,  (\ref{Cuno}) and( \ref {ZSinhECosh})  are the same  as in \cite [p.~136-154] {Edwards:1974cz}).  See  also  \cite [p.~25] {Giovanni Lodone 2021}.

\be \label{Cdue}
C_2(p,\epsilon)=\frac{  i  \psi^{(0)}(p)     }{ 2^5 \pi 3} \left[  15  \left(  1-\frac{ 2 \epsilon}{3}  \right)(1-2 \epsilon) - 5 -9\left( 1-\frac{8 \epsilon}{10}  \right)      \right] - \frac{ i  \psi^{(0)}(p)   (1 +108 \epsilon-12 \epsilon^2    )}{ 2^5 \pi 3}  + 
\ee

$$
 \frac{\psi^{(2)}(p)   5}{2^6 \pi^2}\left(  1-\frac{ 8 \epsilon}{10} \right) -\frac{  \psi^{(2)}(p)  3   }{2^7 \pi^2} \left[      \left(   1-\frac{2 \epsilon}{3} \right)     \left(  1-2 \epsilon \right)    \right] 
 -\frac{\psi^{(2)}(p)   5 }{ 2^7 \pi^2} + \frac{ -i \psi^{(4)} (p) (5-8 \epsilon)  }{2^9 \pi^3 3} + \frac{ i  \psi^{(4)} (p) 5 }{2^9 \pi^3 3}  +\frac{  \psi^{(6)}(p) 5    }{2^7 \pi^4 6!}
$$

\noindent If $\epsilon =0$  the formula is drastically simplified (see  \cite{Edwards:1974cz} pag. 154 bottom).

$$
C_2(p,\epsilon =0) =   \frac{\psi^{(2)}(p)   }{2^6 \pi^2}     +\frac{  \psi^{(6)}(p) 5    }{2^7 \pi^4 6!}
$$

\noindent In \cite [p.~2] {Giovanni Lodone 2021} it  is carried on, even using $M=0$ in (\ref{RDue}),  a significant  numerical  comparison  with :

\be \label {ZMathem}  Z(t)=e^{i\theta(t)}\zeta(1/2+it) \ , \ where \ , \ 
    \theta(t)=   \Im \left[ \ln  \Gamma \left(    \frac{i t}{2} +\frac{1}{4}    \right )      \right]     -\frac{t}{2}  \ln(\pi)    \ ; \ t \rightarrow t-i\epsilon
    \ee

  \noindent see   \cite [p.~119] {Edwards:1974cz} , \cite  {RSiegelWolfram} and appendix \ref{appendix2} .

\noindent This article is structured as follows.

\noindent In section 2  we  recall the Backlund bound   \cite[p.~132]{Edwards:1974cz}  for the  number of zeros of $\xi(s=1/2+\epsilon+it)$ from $t=0$ to $t=T$  : $N(T)$.

\noindent In section 3 we give  a better bound  using   (\ref {ZSinhECosh}) .

\noindent In section 4  %
 we %
  deal about different type of zero configuration.

\noindent  In section 5 we  conclude  summarizing results.%

 \noindent     In appendix   \ref{appendix1}    is computed  the phase of $\Gamma \left( \frac{s}{2}  +1\right) \pi^{-s/2}$   .%

\noindent     In appendix   \ref{appendix2}    is proved the equivalence  of  (\ref {ZSinhECosh}) and  (\ref {ZMathem} ) at $|t|$ big.

\section{Backlund bound  revisited }

\noindent Let us  apply Cauchy theorem  to the closed path of the rectangle defined by:

\be \label {Rect} 
  s_{lower \ left }= -0.5 - i T  \ ;  \   s_{upper  \ right }= 1.5+ iT
  \ee

\noindent    where $T \ne \Im[\rho_k]$ and $\xi(\rho_k)=0$.
For symmetry of $\xi(s)$ with respect to real axis, and, for conjugate symmetry with respect to critical line, we can consider the only path : $ L= L_1+L_2$ where $L_1 : \ (s=1.5) \rightarrow (s=1.5+ i T)$, and, $L_2 : \ (s=1.5+ i T) \rightarrow (s=1/2+ i T)$.  On the whole rectangle total phase variation is $4 \pi N(T)$ because we have $N(T)$ zeros above and below real axis. In $L$ total phase variation is one fourth of  that. So,   on $L=L_1+L_2$,  using  $\angle$ as  phase symbol, we have:
\be \label {CTheorem}  
 \pi N(T) = \angle[ \xi(s) ]_{L_1+L_2} = \angle \left \{ \left[ \frac{\Gamma \left( \frac{s}{2}  +1\right) }{ \pi^{s/2}} \right]  \zeta(s)(s-1) \right\}_{L_1} + \angle \left[ \xi(s)  \right]_{L_2}  
 \ee

\noindent In  \cite[p. ~ 134]{Edwards:1974cz}  it  is stated:

\be \label {Backlund}
\left | N(T) -\frac{T}{2 \pi}\left[ \ln\left ( \frac{T}{2 \pi}\right) -1  \right] - \frac{7}{8} \right |<
B = 0.137 \ln(T) +0.443 \ln(\ln(T))+4.35   \ ; \  \forall T \ge 2
\ee

\noindent  (\ref {Backlund}) is derived by applying Cauchy theorem and Jensen theorem  \cite[p. ~ 40]{Edwards:1974cz}    to rectangle (\ref {Rect}).  See  \cite  [p.~69] {Ingham:1932} )  and  \cite[p. ~133 ]{Edwards:1974cz} .

 \noindent  Jensen's theorem  \cite[p. ~ 40]{Edwards:1974cz}   applies to a rather general kind of  analytic functions,  \cite[p. ~ 133]{Edwards:1974cz}). And this leads 
 to  right main term of (\ref{Backlund}) of the type: $Constant \times \ln(T)$. 
 
 \noindent In (\ref{Backlund}) $Constant= 0.137$. See also \cite  {Hasanalizade2021}.

\section{Refining  Backlund bound  }

\noindent $          LEMMA  \  1:   \ \  $ Suppose $z(x(t))$ is a continuous  complex    function of complex variable   $x(t)$ that describe a curve, not necessarily closed in $t_{in}  \le t \le t_{fin}$.  Besides we suppose {\bf always} $\Im [z(x(t_{in})] =0$. We ask about the  { \bf possible } maximum  absolute total phase change $| \Phi | _{MAX}$,  in the interval
 $t_{in}  \le t \le t_{fin}$ with whatever $\Re [z(x(t)] $%
,  under following  additional  hypothesis:

\begin{itemize}

\item  $\Im [z(x(t)] = 0 \ \forall  t: t_{in}  <  t  \le t_{fin} $ : then   $| \Phi |_{MAX} = \pi$ (if there is a sign change)

\item    $\Im [z(x(t)] \ne 0 \ \forall  t: t_{in}  <  t  \le t_{fin} $ : then   $| \Phi |_{MAX} < \pi$

\item    $\Im [z(x(t)] \ne 0 \ \forall  t: t_{in}  <  t  \le t_{fin} $, but for  only  one $t^*: t_{in}  <  t^*  \le t_{fin} $  we have $\Im [z(x(t^*)] =0  \ $ : 
   then   $| \Phi |_{MAX} < 2 \pi$. %

\end{itemize}
 
\noindent Lemma 1 is self evident if we map $\angle[z(x(t))]$ on the trigonometric circle.

\noindent In order to estimate the phase variation on $L_1$ from (\ref{1p3}), we have: 
 \be \label{FaseSuL1}
 \angle \left[ \xi(s)   \right]_{s=1.5} ^{s=1.5+iT}=
\angle \left[  \Gamma \left( \frac{s}{2}  +1\right)\pi^{-s/2}  \right]_{s=1.5} ^{s=1.5+iT}+\angle \left[(s-1)\zeta(s)    \right]_{s=1.5} ^{s=1.5+iT}
\ee

\noindent If $x \in \Re$ and $x>0$ then  $ \Gamma\left( x \right) \in \Re$,and $ \Gamma\left( x \right)>0$, so
 $\angle\left[ \Gamma\left(\frac{s}{2}+1\right)\pi^{- s/2} \right]_{s=1.5}=   \angle\left[\frac {\Gamma(1.75)} {\pi ^{3/4}} \right] =0$. 
\noindent  So for  (\ref{ImFattDiXiDiSeChi}):

\be \label {FaseSuL1}
\angle\left[ \Gamma\left(\frac{s}{2}+1\right)\pi^{- s/2} \right]_{L_1}  = 
  -\frac {T}{2}+ \frac {T}{2}\ln\left( \frac{T }{2 \pi}  \right)- \frac{\pi}{8}+\frac{\pi}{4}(\epsilon+ 2)  + O(T^{-1}) \  ; \ \epsilon=1
  \ee

\begin{figure}[!htbp]
\begin{center}
\includegraphics[width=1.0\textwidth]{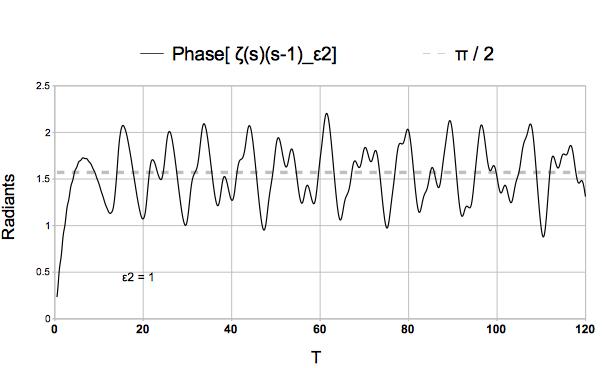} %
\caption{\small {\it   
 $\angle[\zeta(s)(s-1)]_{\epsilon=1}$ in radiants. From figure it is apparent that if   $  \frac{3}{8} \rightarrow \frac{7}{8}$ we can put in (\ref {Backlund2}) $3 \rightarrow 2.5$ at least from $T>10$ as in 
 (\ref {Backlund2_1})  
  }}
\label {Fig1} 
\end{center}
\end{figure}

 \noindent Besides (see fig. \ref {Fig1}):

  \be \label {FaseSuL1DueToZ}
  \left |\angle \left[  \zeta(s)(s-1) \right]_{L_1} \right| < \pi \ \ ;  \and \ : \ \angle \left[  \zeta(s)(s-1) \right]_{s=1.5} =0\ee

  \noindent because from     \cite  [p.~69]  {Ingham:1932} )) :
 \be \label {InghamP69}
  \Re[\zeta(1.5+it)] >1-\frac{1}{2^{1.5}} -\int_2^\infty \frac{du}{u^{1.5}} = 1-\frac{1}{\sqrt{2}}\  \  \ ;  \  \  \ | \Im[\zeta(1.5+it)] | <\frac{1}{\sqrt{2}} \  \  \ ;  \  \  \ \forall t > 0
  \ee
  
 \noindent  So at  least for $T> \frac{1}{2(\sqrt{2}-1)} \approx 1.205$  the $\Im[(s-1)\zeta(s)]_{on \ L_1}=[T \times  \Re[\zeta(s) - 0.5 \times \Im[\zeta(s)]_{on \ L_1} >0 $.
 Direct computation  leads to  a positive value   of $\angle[\zeta(s) (s-1)]$  from $s=1.5$ to $s=1.5+i 1.205$.  See fig  \ref{Fig1}. So  (\ref {FaseSuL1DueToZ}) is true, for Lemma 1 (second point),  on whole $L_1$ .

\noindent In order to estimate  phase difference  on $L_2$ 
 it is possible to use, in (\ref {CTheorem}), at big $T$, instead of properties of a generic analytic function and Jensen theorem, the simple formula  from (\ref {ZSinhECosh}  ):
\be \label {FaseSuL2}
 \angle \left[ \xi(s)   \right]_{s=1.5+ i T} ^{s=1/2+iT}= \left[
 \angle \left[ Z(T,\epsilon)   \right] + \frac {\pi}{4} \epsilon \right]_{\epsilon=1} ^{\epsilon=0}=
 \angle \left[ Z(T,\epsilon)   \right]_{\epsilon=1} ^{\epsilon=0} - \frac {\pi}{4}
 \ee

 \noindent We can simplify $\frac {\pi}{4} \epsilon$ that occurs with opposing sign in (\ref {FaseSuL2}) and in  (\ref {FaseSuL1}), and we can use only  the $\epsilon-$invariant expression:
 
 \be \label {Compensated}
 \left\{
 \angle \left[ \Gamma \left( \frac{s}{2}  +1\right)  \pi^{-s/2} \right] -\frac {\pi}{4} \epsilon
 \right\}_{\epsilon-inv.}= \frac{T}{2} \ln \left ( \frac{T}{2 \pi e} \right)  +\frac{3}{8}\pi   +O(T^{-1}) \ ; \ \forall \epsilon 
 \ee
 
 \noindent  for (\ref {FaseSuL1}), and only  $ \angle \left[ 
 Z(T,\epsilon)    \right]_{\epsilon=1} ^{\epsilon=0}  $ for   \ref{FaseSuL2}, like in fig. \ref  {Fig6}  .

\be \label {costr1}  LEMMA  \ 2 : \ \ \ \ \ \
\left|  \angle \left[ Z(T,\epsilon)   \right]_{L_2} \right| <2 \pi
\ee

\noindent  Proof. We focus on  $\Im \left[  Z(T,\epsilon) \right]_{L_2}$.  At constant $T$  we have a sum of decreasing  hyperbolic sines  as $\sum_{decr}A_j \sinh(a_j \epsilon)$, and  increasing  hyperbolic sines  as $\sum_{incr}B_j \sinh(b_j \epsilon)$, both are null for $\epsilon=0$,   plus  terms that, at first,  can be neglected. %

\noindent Coefficient of    $\sum_{decr}A_j \sinh(a_j \epsilon)$ and $\sum_{incr}B_j \sinh(b_j \epsilon)$ are from (\ref{ZSinhECosh}) :

\be \label {aAbB}
 |a_n| , |b_n| =\left|     \ln\left( \sqrt{\frac {T}{2 \pi n^2}  }\right)   \right| \  ; 
  \  |A_n| , |B_n| =   \left|  
 \frac{  1  }{\sqrt{n}}
\sin \left(   T  \  \ln\left(    \sqrt{\frac {T}{2  e \pi n^2} } \right) - \frac{\pi}{8}  \right )   
 \right|
 \ee

\noindent Consequently in $L_2$, for $\Im \left[  Z(T,\epsilon) \right]_{L_2}$ ( see (\ref{ZSinhECosh}) ),  it may happen to have  no zero at all, or, just a  single zero (excluding end  points). No other possibility are left as  $- \sum_{decr}A_j \sinh(a_j\epsilon)$, and  $\sum_{incr}B_j \sinh(b_j \epsilon)$   are zero at $\epsilon=0$ and are  both increasing for $\epsilon >0$.

\noindent%
Besides  $ C_{k}(p,\epsilon)$  or  $ |Err(t)|$ do not influence above conclusion.

\begin{figure}[!htbp]
\begin{center}
\includegraphics[width=1.0\textwidth]{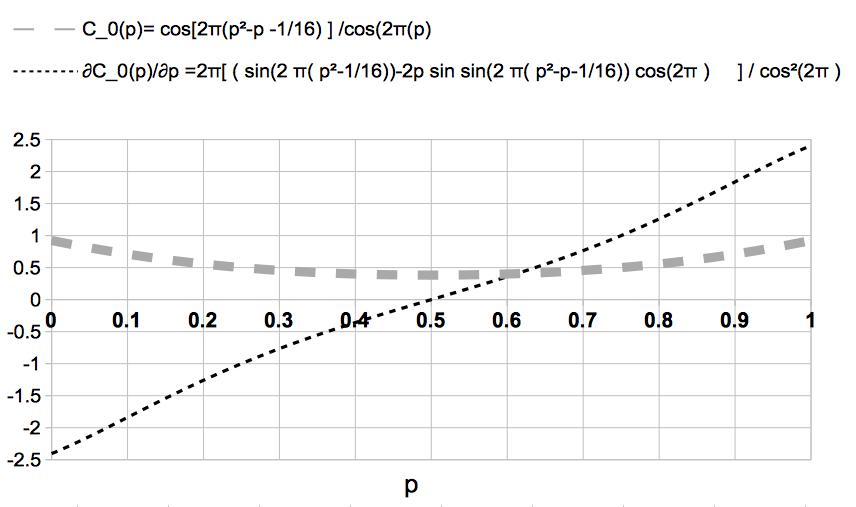}.
\caption{\small {\it   
 $\Im[R_n(t,\epsilon) ]$ is sum of polynomials in $\epsilon$, without $\epsilon^0$ term,  with coefficients beeing products of constants and factors like
$\Psi(p)^{(n)}$.  Where  
$\Psi(p)^{(n)}$if $n$ even is  similar in shape  to $\Psi(p)^{(0)} = C_0(t)$, if $n$ odd it is  similar in shape  to $\Psi(p)^{(1)} = \frac{\partial \Psi(p)^{(0)} } { \partial p} = \frac{4 \pi}{\epsilon}  \Im[C_1(t, \epsilon)]$. See (\ref {RDue}  ), (\ref {Cuno}  ) , and,  (\ref {Cdue}  ).
  }}
\label {Fig8} 
\end{center}
\end{figure}  

\noindent From (\ref{ZSinhECosh})  %
we find   $\Im[R(t,\epsilon) 
 ]$. %
   Following  
  (\ref {RDue}  ), (\ref {Cuno}  ) , and,  (\ref {Cdue}  ) we have  (see also fig. \ref {Fig8} ):

\noindent  For $n=1$ we have $ \Im[C_1(t, \epsilon)] = -  \frac{\epsilon} {4 \pi}  \frac{\partial \Psi(p)^{(0)} } { \partial p}$. 

\noindent For $n=2$ we have:

  $$ \Im[C_2(t, \epsilon)] =      \psi^{(0)}(p)   \left[   \frac{ 1    }{ 2^5 \pi 3} \left[  15  \left(  1-\frac{ 2 \epsilon}{3}  \right)(1-2 \epsilon) - 5 -9\left( 1-\frac{8 \epsilon}{10}  \right)      \right] - \frac{   (1 +108 \epsilon-12 \epsilon^2    )}{ 2^5 \pi 3}   \right]+      \psi^{(4)} (p)   \frac{8 \epsilon}{2^9 \pi^3 3}     $$

\noindent Afterword there is the $(-1)^{N-1}\left(\frac{2 \pi}{t}  \right)^{1/4+n/2} $ factor, see (\ref{RDue}),  that decreases with increasing  $n$. 
We are in $0 \le \epsilon \le 1$ interval, so the $\epsilon^n$ powers ($n \ge 1$), %
behave like the increasing or decreasing $\sinh( )$ functions above, and so they  can join the above sums of increasing or decreasing functions that vanish for $\epsilon=0$, so,  confirming previous result.

 \noindent $\Im[Err(t,\epsilon)] <e^{-t/10}$ for $ t $ big  is of course very small  (see (\ref{ZSinhECosh})).  Besides his variation with $\epsilon$ is very slow  as it is connected with integrals relatively far from critical line. See \cite [p.~28-38]  {Giovanni Lodone 2021}. So from Lemma 1 (third  point ), the inequality  (\ref {costr1}), i.e. Lemma 2, is true.

\noindent Thanks to (\ref{FaseSuL1}), (\ref {FaseSuL1DueToZ}), 
and (\ref {costr1}),  we have (neglecting $O(T^{-1}) $ for big $T$, see  also \cite  [p.~69]  {Ingham:1932}, and appendix \ref {appendix1}):
$$
\angle \left [ \xi(s) ]_{L_1+L_2}\right]= \left[ \frac{T}{2} \ln \left ( \frac{T}{2 \pi e} \right)  +\frac{3}{8}\pi  \right] +\angle \left[  \zeta(s)(s-1) \right]_{L_1} + \angle \left[  Z(T,\epsilon) \right]_{L_2} 
$$
$$ <\left[ \frac{T}{2} \ln \left ( \frac{T}{2 \pi e} \right)  +\frac{3}{8}\pi  \right] + \pi+2\pi$$

 \noindent  So   Backlund estimate (\ref {Backlund}) can be improved (see fig. \ref {Fig10}  )  at least with: 

\be \label {Backlund2}  THEOREM  \ : \
\left |
 N(T) -\frac{T}{2 \pi}    \ln \left ( \frac{T}{2 \pi e} \right) - \frac{3}{8}
  \right |
 < 3 \ ;   \ \ \forall T >10, and , \ne \Im[\rho_k] \ ; \ \xi(\rho_k)=0
\ee

 \noindent  Proof:   phase variations  of $\left| \angle \left[  \zeta(s)(s-1) \right]_{L_1} + \angle \left[  Z(T,\epsilon \right]_{L_2}  \right|$  in  (\ref{CTheorem}  )   are  surely $< 3\pi$. 

  \noindent In  (\ref {Backlund2}) $  \frac{3}{8}$ replaces $\frac{7}{8}$ in (\ref {Backlund})  because we have gathered  $\zeta(s)(s-1)$ instead of  leaving $\zeta(s)$ alone.  See fig. \ref  {Fig1} .  So we can write:
  
\be \label {Backlund2_1}
\left |
 N(T) -\frac{T}{2 \pi}    \ln \left ( \frac{T}{2 \pi e} \right) - \frac{7}{8}
  \right |
 < 2.5 \ ;   \ \ \forall T >10, and ,  \ne \Im[\rho_k] \ ; \ \xi(\rho_k)=0
\ee

   \noindent Notice that  {\bf   this bound  is no $T-$dependent  }  like  (\ref{Backlund}).  %
 See  fig \ref {Fig10} .  
 Direct verification till, at least,  $2\times 10^6$ first zeros (\cite  {Odlyzko})  complies with (\ref{Backlund2}). 


  %

\section{%
Analysis of zero configurations  by (\ref {ZSinhECosh}). 
}
\label {Cases}
If $N_{off}(T)$ are the hypothetical  off-critical line zeros till $T$, while  $N_{on}(T)$ are the on-critical line zeros :  $N(T)=N_{on}(T) +N_{off}(T)$, but along $L$  :

\be \label {RevisedBacklundBound2_}
\left| N_{on}(T)  + 2 N_{off}(T)  -
 \left \{ \frac {T}{2\pi} \ln\left( \frac{T }{2 e \pi}   \right) + \frac{3}{8} \right\} \right| < 3
 \ee
Notice that, for symmetry properties of $\xi(s)$  (i.e. $\xi(1-s)=\xi(s)$, and,  $\tilde{\xi(s)}=\xi(\tilde{s}) $), an unitary increase of $N_{off}(T)$   means a $2\pi$ increase of the total phase on $L$ path. While an unitary increase of   $N_{on}(T)$   means  only a $\pi$ increase of total phase on $L$ path.
 This is apparent  as two  internal  $\frac {1}{2}  \pm it_0 $  zeros  give  $4 \pi$ total phase increment on (\ref {Rect}), while four internal  $\frac {1}{2} \pm \epsilon_0  \pm it_0 $  zeros give  $8 \pi$ total phase increment on same (\ref {Rect}), and, phase on $L$ is $1 /4$  of total phase on  (\ref {Rect}).

\begin{figure}[!htbp]
\begin{center}
\includegraphics[width=1.0\textwidth]{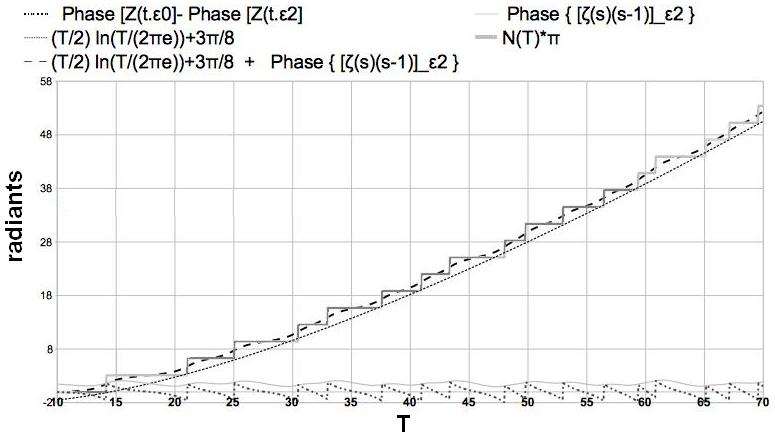} %
\caption{\small {\it   
 Comparison among  $\pi N(T)=
\angle\left[ \Gamma\left(\frac{s}{2}+1\right)\pi^{- s/2} \right]_{L_1}+\angle \left[  \zeta(s)(s-1) \right]_{L_1} +  \angle \left[  Z(T,\epsilon) \right]_{L_2} $  and their  single constituent $\angle\left[ \Gamma\left(\frac{s}{2} \ ; \ 1\right)\pi^{- s/2} \right]_{L_1} \ $ by (\ref  {Compensated} ) $ ;  \ \angle \left[  \zeta(s)(s-1) \right]_{L_1} \ ; \ \angle \left[  Z(T,\epsilon) \right]_{L_2} $. The sum of the three phase differences gives exactly  the number of zeros  at $T$ multiplied by $\pi$. $\epsilon_0=0$, and,  $\epsilon_2=1$ .  In this, and ,following figures,
(\ref{ZSinhECosh}) is computed with  $M=0$ in (\ref{RDue}).
}}
\label {Fig6} 
\end{center}
\end{figure}

\begin{figure}[!htbp]
\begin{center}
\includegraphics[width=1.0\textwidth]{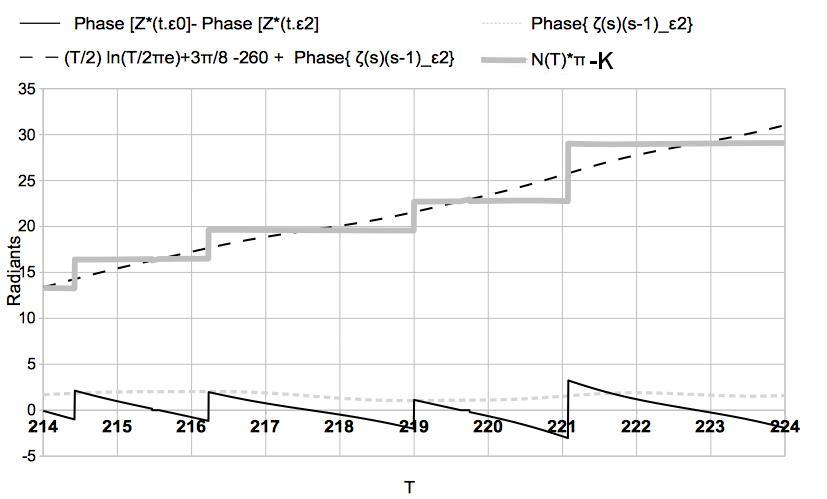} %
\caption{\small {\it   
First off critical line zero found with $Z^*(T,\epsilon)$ computed as  (\ref{ZSinhECosh})  with $R(t,\epsilon)=0$. The step at $T\approx 221$ is of $2\pi$. The perturbation on $\pi N(T) -K$ ( $K=260$), with respect to $Z(t,\epsilon)$ case is shown in fig. \ref  {Fig9} .
  }}
\label {Fig5} 
\end{center}
\end{figure}

\begin{figure}[!htbp]
\begin{center}
\includegraphics[width=1.0\textwidth]{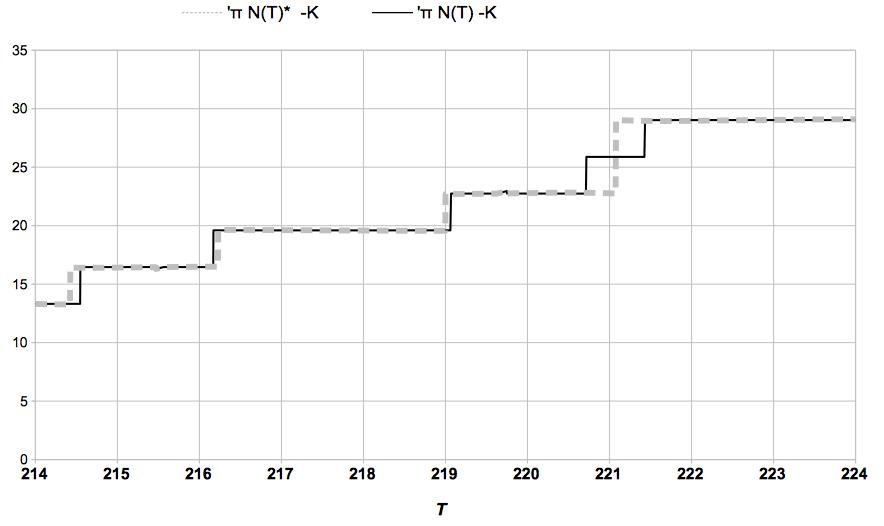} %
\caption{\small {\it   
 Zero count comparison  between  $Z(T,\epsilon)$ and his perturbed version $Z^*(T,\epsilon)$. In  the perturbation $Z\rightarrow Z^*$  expression (\ref {RDue}), in  (\ref{ZSinhECosh}) , has been erased, i.e.  we have arbitrarily forced   $R(t,\epsilon)=0$. 
Notice that two on-critical line zeros in  $\pi N(T)$,  around $T=221$, are transformed in two off-critical line producing a  $2 \pi$ jump in  $\pi N^*(T)$.%
But the zeros count is almost unaffected by the perturbation $Z\rightarrow Z^*$ far from  $T=221$. The constant $K =  260$.
  }}
\label {Fig9} 
\end{center}
\end{figure}

\noindent Nowadays we know only $N_{on}(T)$ type  zeros   \cite {Gourdon:2004cz}. See fig. \ref{Fig6}. But if we force $R(t,\epsilon)=0$ in (\ref{ZSinhECosh}),  it  produces   a function,  here called $Z^*(t,\epsilon)$, very similar to $Z(t,\epsilon)$
, that, 
gives rise to some  off-critical line zeros.  The comparison between  $Z^*(t,\epsilon)$, and,  $Z(t,\epsilon)$ can be carried on only locally  to avoid discontinuities of $Z^*(t,\epsilon)$ at $T=2 \pi N^2$.  In  figg. \ref {Fig5} , and \ref {Fig9}  the first of these fictitious zeros arises for $T\approx 221$. Many others follow at higher $T$. Perhaps it is worth to mention that first Lehmer point  ($T\approx 7005.2$  \cite[p.~178] {Edwards:1974cz})  ) happens to be off critical line, for $Z^*(t,\epsilon)$ of course. But, strangely, not second one ($T\approx 17143.8$  \cite[p.~179] {Edwards:1974cz})  ).

\noindent In any case the phenomenology of an hypothetical off-critical line zero for $Z(t,\epsilon)$ or $\xi(s)$  (if it exists), seen by the  $\Re[\xi(s)]=0$, and,  $\Im[\xi(s)]=0$ curves, may  be similar  to fig. \ref {Fig2} based on $Z^*(t,\epsilon)$. Compare with fig. \ref  {Fig3} from $Z(t,\epsilon)$ in same interval.

\noindent Observing figg. \ref  {Fig10} or \ref    {Fig4}, one can wonder  if bound  2.5 in (\ref {Backlund2_1}) could be reduced, or, if    an interval for the minimum B value  can be foreseen. 
From computation like in fig. (\ref {Fig7} ) applied around $T\approx 7005.2$ it is easily seen  that the value 1 in place of 2.5 is  inadequate. So Backlund bound B from :
$$ B =  0.137 \ln(T) +0.443 \ln(\ln(T))+4.35 \ \ ; \  \forall T \ge 2
$$

becomes a T-invariant number B with at least following constraints: 
$$ 1 < B \le 2.5 \ \ ; \  \forall T > 10 $$

\begin{figure}[!htbp]
\begin{center}
\includegraphics[width=1.0\textwidth]{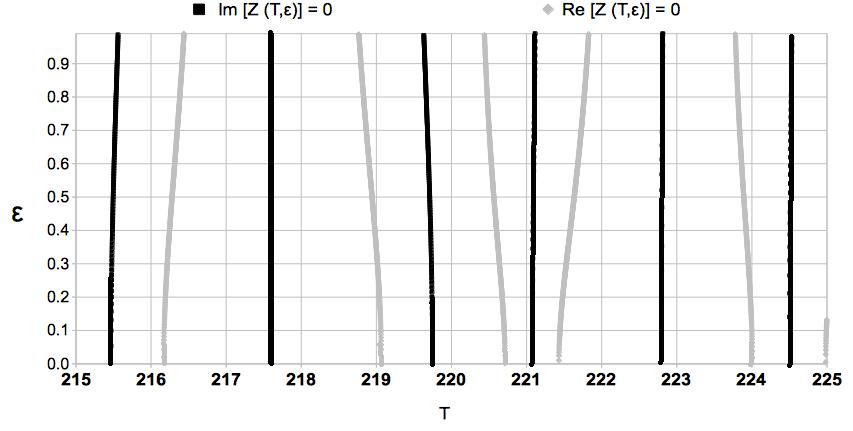} %
\caption{\small {\it   
Locus of points with  $\Re[Z(T,\epsilon)]=0$ or  $\Im[Z(T,\epsilon)]=0$.
  }}
\label {d} 
\end{center}
\end{figure}

\begin{figure}[!htbp]
\begin{center}
\includegraphics[width=1.0\textwidth]{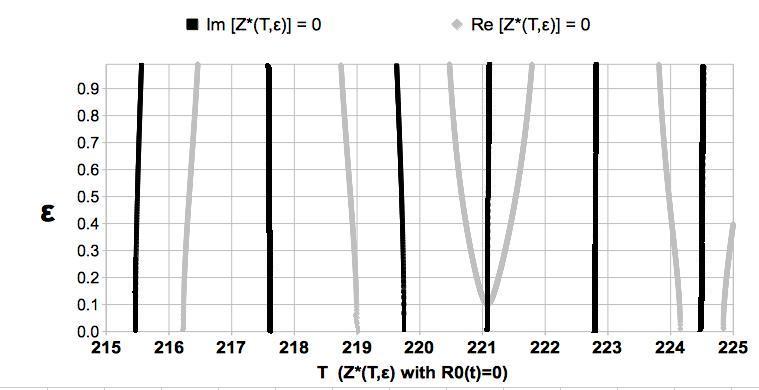} %
\caption{\small {\it   
Locus of points with  $\Re[Z^*(T,\epsilon)]=0$ or  $\Im[Z^*(T,\epsilon)]=0$. Notice the off critical line zero at $T\approx 221$.
  }}
\label {Fig2} 
\end{center}
\end{figure}

\begin{figure}[!htbp]
\begin{center}
\includegraphics[width=1.0\textwidth]{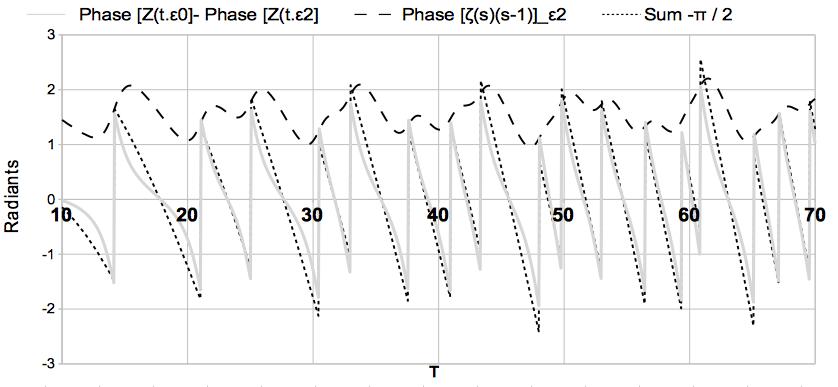} %
\caption{\small {\it   
Magnified view of $\angle[\zeta(s)(s-1)]_{\epsilon=1}$ and $\angle[Z(t,\epsilon)]_{L_2}  $  in fig. \ref {Fig6}. Notice that the sum 
 $\angle[\zeta(s)(s-1)]_{\epsilon=1}+\angle[Z(t,\epsilon)]_{L_2} -\frac{\pi}{2} $ is  rather linear. It is connected to the flat part of $\pi N(T)$ in fig. \ref {Fig6}. For (\ref{Backlund2_1} ) this sum is  always in the strip $\pm2.5 \pi$. Cases analysis suggests it cannot be narrower than $\pm \pi$.The $\pi$ vertical jumps are due to $\xi(s)$  zeros. 
  $\epsilon_0=0$, and,  $\epsilon_2=1$ . }}
\label {Fig7} 
\end{center}
\end{figure}

\begin{figure}[!htbp]
\begin{center}
\includegraphics[width=1.0\textwidth]{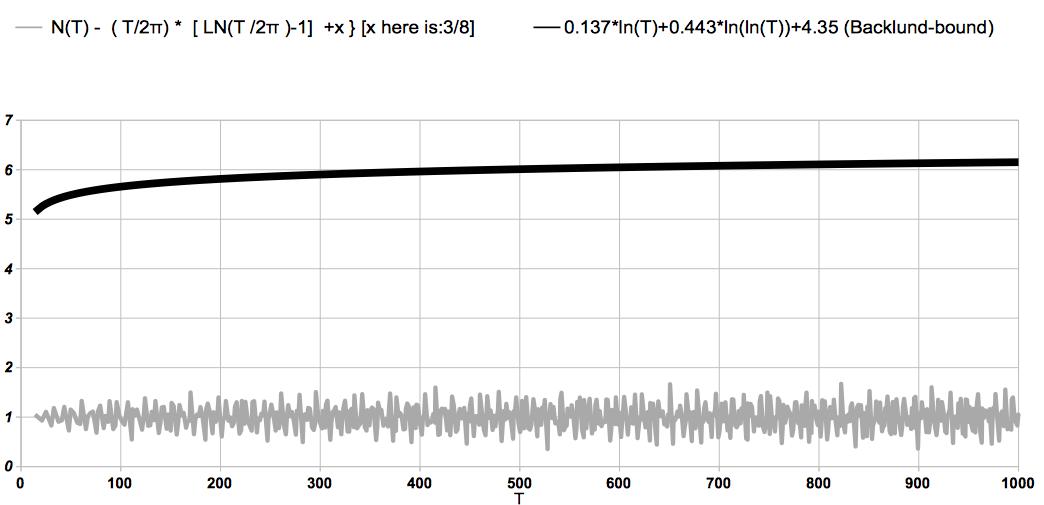} %
\caption{\small {\it   
 Also at higher t values, from an experimental point of view, Backlund bound seems too pessimistic. 
. Expression : $N(T) -\frac{T}{2 \pi}\left[ \ln\left ( \frac{T}{2 \pi}\right) -1  \right] - \frac{3}{8} $, here evaluated at $ T  = \Im[\rho_k] \ ; \ \xi(\rho_k)=0$. This sampling captures only upper points of, for example,  fig  \ref {Fig4}. 
Data on $N(T)$ from \cite {Odlyzko}. }}
\label {Fig10} 
\end{center}
\end{figure}

\begin{figure}[!htbp]
\begin{center}
\includegraphics[width=1.0\textwidth]{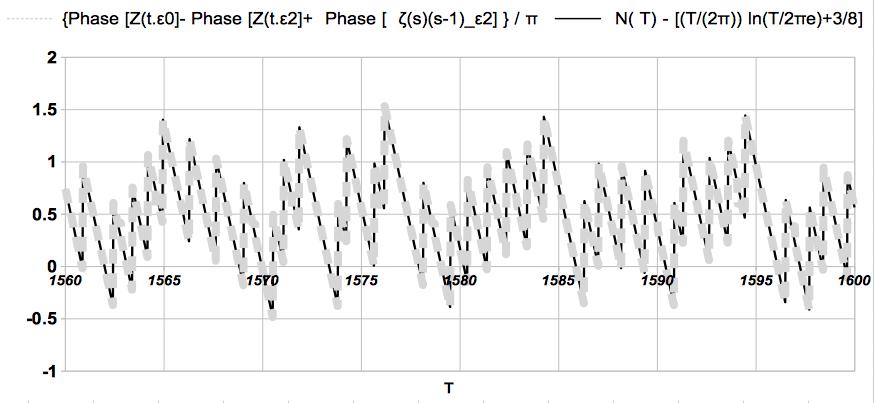} %
\caption{\small {\it   
 %
. Expression :$N(T) -\frac{T}{2 \pi} \ln\left ( \frac{T}{2 e \pi}\right)    - \frac{3}{8} $ is  here   resampled at $\Delta T = 0.027$, and compared  with $\frac{\angle \left[  \zeta(s)(s-1) \right]_{L_1} + \angle \left[  Z(T,\epsilon) \right]_{L_2}}{\pi}$ : $\epsilon_0=0$, and,  $\epsilon_2=1.5$ .
Data on $N(T)$ from \cite {Odlyzko}. }}
\label {Fig4} 
\end{center}
\end{figure}

\section { Conclusions}

\noindent The use of (\ref{ZSinhECosh}) representation of  $\angle[\xi(s)[$ allows us to  revise the Backlund bound from (\ref   {Backlund}  )  to  (\ref {Backlund2_1}); the latter  explains better  fig. \ref  {Fig10}  and similar tests at higher $t$ values.

\noindent Besides (\ref{ZSinhECosh}) has been used  in section   \ref {Cases}         for analysing  possible  types of zeros of  $Z(t,\epsilon)$, an object closely linked to $\xi(s)$. See also appendix  \ref  {appendix2}  .

\section*{Acknowledgments}

\noindent A thankful thought to dear professor Richard Bruce Paris, unfortunately  departed about two years ago, for his kind support and  precious suggestions for  \cite  {Giovanni Lodone 2021}  which  is the  basis of this paper. %
I also thank Paolo Lodone for  our useful discussions and for contributing in some points of this work.


\vspace{0.3cm}

\appendix

\section{ %
Phase of $\Gamma \left( \frac{s}{2}  +1\right) \pi^{-s/2}$  %
 } \label{appendix1}

\noindent From (\ref {1p3}):

\be \label {FattorePerXIAChi}
 \angle[\xi(s)] =  
  \angle \left[  \left( \pi   \right)^{-\frac{s}{2}} \Gamma  \left( \frac{s+2}{2} \right)  \right] \ + 
  \angle  \left[  (s-1)  \zeta(s) \right]
\ee
%
\noindent If we consider the factor $(s-1)$ attached to $\zeta(s)$ 
 \noindent we can use Stirling formula ( \cite[p.~109,112]{Edwards:1974cz})  putting : 
$z+1=  \frac{s+2}{2} \rightarrow z=\frac{2\epsilon+4-3}{4}+i \ \frac{t}{2}$,  so we have:

\be\label {Pi&Gamma3}
\ln( \Gamma(z+1) )  +\ln\left[\left( \frac{\pi}{q} \right)^{ - \frac{s}{2}} \right]=
\ln( \Gamma(z+1) )  -\frac{s}{2} \ \ln\left(\pi\right) =
\ee
$$\ln\left( e^{-z}z^{   z+\frac{1}{2}   }(2\pi)^{\frac{1}{2}}  \right )  +  \left( \sum_{k=1}^{K-1} \frac{B_{2k}   }{  2 k (2k-1) z^{2k-1}} \right) +R_{2K}(z)  -\ln\left(   \pi    \right)\frac{s}{2} =\Re_1+\Re_2+\Re_3 + i(\Im_1+\Im_2+\Im_3)
$$
Let us call    $\Im_1$    the imaginary  of  asymptotic  (for $t $ big) part of   $\ln\left( e^{-z}z^{   z+\frac{1}{2}   }(2\pi)^{\frac{1}{2}}  \right ) $. And $\Re_1$ is the real part.  Let us call    $\Im_2$   the left imaginary part  of  $\ln\left( e^{-z}z^{   z+\frac{1}{2}   }(2\pi)^{\frac{1}{2}}  \right ) $ wich goes to zero  for $t\rightarrow \infty$. And $\Re_2$ is the real part.And finally let us call $\Re_3+i \Im_3 =    \left( \sum_{k=1}^{K-1} \frac{B_{2k}   }{  2 k (2k-1) z^{2k-1}} \right) +R_{2K}(z)   $.

\noindent The $B_j$  are the Bernoulli numbers  ( \cite[p.~11] {Edwards:1974cz} ), that vanish for odd $j$ %

The modulus of  the error term $|R_{2K}(z)|$ is bounded by: 
\be \label  {RestoInApp}
 |R_{2K}(z)| < \left(    \frac{B_{2K}   }{  2 K (2k-1) z^{2K-1}}          \right)
\frac{1}{       \left[     \cos\left ( \frac{\arg(z)}{2}\right)         \right]^{2K}        }
\ee
where $\arg(z)$ is taken in the interval: $-\pi <\arg(z) <\pi$ (see  \cite[p.~112]{Edwards:1974cz}.

 \noindent  In (\ref{Pi&Gamma3}) we may put : $z=\frac{2\epsilon+1}{4}+i \ \frac{t}{2}$ , so,   leaving apart only $\Im_3$ for now, we have:

 $$
\ln\left[   \left( \frac{\pi}{q}   \right)^{-\frac{s+2_1}{2}} \Gamma  \left( \frac{s+2}{2} \right) \right] 
\approx  \Re_1 + \Re_2+ i ( \Im_1 + \Im_2) = %
 $$

\noindent   Taking into account that:
 $
  \arctan \left[ \Im=\frac{ t}{2} ; \Re = \frac {\epsilon+ 2}{2}-\frac{3}{4} \right]=\frac{\pi}{2}+\arctan \left[  \frac{3-2(\epsilon+2)}{2t}  \right] 
 $, we get:
  \be \label {NoApproxSalvo Bernoulli2}
 \frac{3}{4} - \frac {\epsilon+ 2}{2}-\frac{i t}{2} +\left (\frac {\epsilon+ 2}{2}- \frac{1}{4} + \frac{i t}{2}\right) \left\{     \ln  \sqrt{\left( \frac {\epsilon+ 2}{2}-\frac{3}{4} \right)^2+\frac{t^2}{4}}  +i \left \{ \frac{\pi}{2}+\arctan \left[  \frac{3-2(\epsilon+2)} {2t}  \right]   \right \}   %
  \right\}
  +
  \ee

$$
 \frac{\ln(2 \pi)}{2}-\ln\left(\pi \right)\left(    \frac {\epsilon}{2} +\frac{1}{4}   \right)  -i\frac{t}{2}\ln\left(  \pi  \right)
$$

\noindent   From which we find the  the asymptotic  imaginary parts  $\Im_1 $ :

 \be \label {ImFattDiXiDiSeChi}
\Im_1
 = 
  -\frac {t}{2}+ \frac {t}{2}\ln\left( \frac{t }{2 \pi}  \right)- \frac{\pi}{8}+\frac{\pi}{4}(\epsilon+2)
 \ee

\noindent  for $t$ big: 
  $$ \Im_2  \approx \frac{t}{4}  \left( \frac{2(\epsilon+ 2)-3}{2t}  \right)^2 +
   \left(    \frac {\epsilon+ 2}{2} -  \frac{1}{4}   \right)  \left(  \frac{3-2(\epsilon+2)}{2t}      \right) =O(t^{-1})
  $$

 While $\Im_3=  \Im \left[     \left( \sum_{k=1}^{K-1} \frac{B_{2k}   }{  2 k (2k-1) z^{2k-1}} \right) +R_{2K}(z)  \right]$ for $K=3$.So $k = 1 , and \ 2$

\be \label {Im3TillK2}
\Im_3 -R_{2K}(z)= \frac {-1}{6 t\left[ 1+\left( \frac{2 \epsilon+2 2-3 }{2t}  \right) ^2\right]}
-\frac {1}{45 t^3 \left[1 +\left(\frac{2\epsilon+22-3}{2t}  \right)^2 \right]^3}
+
\frac {(2\epsilon+22-3)^2}{60 t^5 \left[1 +\left(\frac{2\epsilon+22-3}{2t}  \right)^2 \right]^3}   =  O(t^{-1})  %
\ee

\section{ %
About $Z(t-i\epsilon)$ ( \ref  {ZMathem}  ) and $Z(t,\epsilon)$  (\ref  {ZSinhECosh} )
 } \label {appendix2}

 \noindent  Formula (\ref {ZMathem}), according to  (\cite  {RSiegelWolfram}),  { \it `' is an analytic function of t except for branch cuts on the imaginary axis running from $\pm \frac{i}{2}$  to $\pm i \infty $''}. 
 
 \noindent Z(t)  of (\ref {ZMathem}) is implemented 
  in  {\it Wolfram Mathematica}  by RiemannSiegelZ function.  So  $Z(t-i\epsilon)$ can be compared with $Z(t,\epsilon)$ in  (\ref {ZSinhECosh}  ) which is analytic too at high $t$ and limited $\epsilon$.

\noindent  Let us see it  (for $|t|$ so big that  $t^{1/4} \approx [s(s-1)]^{1/8} $) :

  $$F(t) e^{i \frac{\pi}{4}  \epsilon}   = \left(\pi/2 \right)^{0.25}t^{\frac{7}{4}} e^{- \frac{\pi}{4}  t} e^{i \frac{\pi}{4}  \epsilon}  =_{t \ big}    \left(\pi/2 \right)^{0.25}[s(s-1)]^{\frac{7}{8}} e^{ i \frac{\pi}{4} ( s-1/2) }$$
  
  \noindent So $Z(t,\epsilon)$ in  (\ref {ZSinhECosh}  )  is analytic  because ratio of analytic functions ( $|t|>>1 \ ; 0\le \epsilon <2$).
 \noindent Below is  explained  why numerical evaluation of  $Z(t,\epsilon)$ and $ Z(t-i\epsilon)$ gives almost identical result.

\noindent  It happens to meet  often expressions like  $\Gamma\left(\frac{s}{2}\right) $, for example in (\ref {ZMathem} ), or  in (\ref{1p3}). Or expressions like $\Gamma\left(\frac{s}{2}+1\right) $ or $\Gamma\left(\frac{(1-s)}{2}+1\right) $ like in 
 (see \cite[p.~138]{Edwards:1974cz}
 :
\be \label{StartingP}
 -\xi(s) = (1-s)\Gamma \left(  \frac{s}{2} +1\right)\pi^{\frac{-s}{2}} \left( \sum_{n=1}^{n=N}n ^{-s}  \right) +
(+s)\Gamma\left(  \frac{1-s}{2} +1\right)\pi^{-\frac{1-s}{2}} \left( \sum_{n=1}^{n=N}n ^{-(1-s)}  \right)  +  
\ee

$$+ \frac{ (+s)\Gamma\left(  \frac{1-s}{2} +1 \right)\pi^{\frac{-(1-s)}{2}}   }{  (2\pi)^{s-1} \,  2 \sin(\pi s/2) \, 2\pi i   }\int_{C_{N}} \frac{(-x)^{s-1}e^{-Nx}dx}{e^x-1} 
$$
wich is the starting point for     \cite  {Giovanni Lodone 2021}.
                        \noindent It is then convenient to gather all this computations.

\begin{table}[bht]
\begin{center}
\begin{tabular}{c|c|c|c|c|c|cl}
$Cases$     & $ \ln[\Gamma( . . .)]$    & $\alpha$     & $\beta   $     &     $\pi \ $ power     &$ \angle [s_{factor}  ]$     \\
\hline
$1$ & $\ln\left[\Gamma\left(\frac{s}{2}+1\right) \right]$ & $1$& $1$ &
$\left(-\frac{1}{4}-\frac{\beta(\epsilon+i t)}{2} \right) \ln(\pi)$  
& $\angle[(1-s)]_{t \ big} = -\frac{\pi}{2}$   \\
\hline
$2$ & $\ln\left[\Gamma\left(\frac{1-s}{2}+1\right) \right]$ &$1$ & $-1$ &  
$\left(-\frac{1}{4}-\frac{\beta(\epsilon+i t)}{2} \right) \ln(\pi)$ & $\angle[(+s)]_{t \ big} = +\frac{\pi}{2}$ \\
\hline
$3$ &$\ln\left[\Gamma\left(\frac{s}{2}\right) \right]$ & $-3$& $1$ &
$\left(-\frac{1}{4}-\frac{\beta(\epsilon+i t)}{2} \right) \ln(\pi)$ 
& $\angle[(+s)(1-s)]_{t \ big} = 0$ \\
\end{tabular}
\end{center}
\caption{\small  With parameters in table the canonical Stirling formula $\ln( \Gamma(z+1) ) = -(z) +\left(z+\frac{1}{2}\right) \ln(z) +0.5 \ln(2 \pi) + O(z^{-1} )$ has same espression in above three cases, %
if $z= \frac{\alpha}{4}+\frac{\beta \epsilon}{2}+\frac{\beta i t}{2} $. Also $\pi-$powers, i.e. $\pi^{-s/2}$ or $\pi^{-(1-s)/2}$,  assume same expression if we use  appropriate $\alpha,\beta$ parameters.
}
\label {3Cases}
\end{table}

 \noindent Applying Stirling formula ( and neglecting for the moment $O(z^{-1} )$) : 
   $$\ln \left[
 \Gamma
 \left (
 \frac{\alpha}{4}+\frac{\beta \epsilon}{2}+\frac{\beta i t}{2} +1
 \right)\right]= 
 -\left (
 \frac{\alpha}{4}+\frac{\beta \epsilon}{2}+\frac{\beta i t}{2} 
 \right) + \ln(\sqrt{2 \pi})
 $$
 $$
 +\left (
 \frac{\alpha}{4}+\frac{1}{2}+\frac{\beta \epsilon}{2}+\frac{\beta i t}{2} 
 \right)
  \left [ \ln\sqrt{\left (
 \frac{\alpha}{4}+\frac{\beta \epsilon}{2} 
 \right)^2 +\left (
 \frac{\beta  t}{2} 
 \right)^2  } +i \left(\frac{\beta \pi}{2}     
      -\beta \arctan \left(    \frac{(\beta \epsilon)/2+  \alpha/4}{(\beta t)/2} \right)    \right) \right]
 $$ 
  Notice that  $-\left (
 \frac{\alpha}{4}+\frac{\beta \epsilon}{2}
 \right)$ and $\frac{\beta i t}{2}  \times ( - i  )
 \beta \arctan \left(    \frac{  (\beta \epsilon)/2+  \alpha/4  }{  (\beta t)/2} \right)    $ simplifies $\forall \beta$  for  $t \rightarrow \infty$, so we get :
 $$- \frac{\beta i t}{2} +\left (
 \frac{\alpha}{4}+\frac{1}{2}+\frac{\beta \epsilon}{2} 
 \right) \ln\left( \frac{t}{2}\right) +\frac{\beta i t}{2} \ln\left( \frac{t}{2}\right) +
 \frac{\beta i \pi}{2} 
 \left (
 \frac{\alpha}{4}+\frac{1}{2}+\frac{\beta \epsilon}{2} 
 \right) - \frac{\beta^2 \pi}{4} t + \ln(\sqrt{2 \pi})
 $$

\noindent Then we conclude:
   
\be \label {RealPart}
\Re \ln \left[
 \Gamma
 \left (
 \frac{\alpha}{4}+\frac{\beta \epsilon}{2}+\frac{\beta i t}{2} +1
 \right)\right]=
 \left[  \frac{\alpha}{4} + \frac{1}{2} +  \frac{\beta \epsilon}{2}\right] \ln\left(\frac{t}{2} \right)
 +\ln(\sqrt{2 \pi  })-\frac{\beta^2\pi}{4}t + O(t^{-1} )
\ee

\be \label {ImmPart}
\Im \ln \left[
 \Gamma
 \left (
 \frac{\alpha}{4}+\frac{\beta \epsilon}{2}+\frac{\beta i t}{2} +1
 \right)\right]=
 \beta\left[\frac{t}{2} \ \ln\left( \frac{t}{2 e}  \right) + \frac{\pi(\alpha +2)}{8} \right]
 + \frac{ \beta^2  \pi}{4} \epsilon + O(t^{-1} )
\ee
   
 \noindent From  now on  we choose %
  $ \beta=1 , \alpha=-3$ for $\Gamma\left(\frac{s}{2}\right)$ in    (\ref {ZMathem} ), and  (\ref{1p3}) :
 
 \be \label {Zt2}
 Z(t) = e^{ i \left(   \beta\left[\frac{t}{2} \ \ln\left( \frac{t}{2 e}  \right) + \frac{\pi(\alpha +2)}{8} \right]
 -\frac{\beta(i t)}{2} \ln(\pi) + O(t^{-1} ) \right) } \ 
 \frac { (- \xi(1/2 +it))} 
 {  
 \Gamma \left(   \frac{ 1/2+it} {2}      \right) \pi^{\left(-\frac{1}{4}-\frac{\beta(i t)}{2} \right)} 
  \frac { (1/2+it) (1/2-it)} {2} 
   }
 \ee
 
 When we put  $t \rightarrow t- i\epsilon$,  we have for the ratio in (\ref{Zt2}) :
 $$
  \frac { (- \xi(1/2 +it+\epsilon))} 
 {  
 \Gamma \left(   \frac{ 1/2+it +\epsilon} {2}      \right) \pi^{\left(-\frac{1}{4}-\frac{\beta(\epsilon+i t)}{2} \right)} 
  \frac { (1/2+it + \epsilon) (1/2-it -\epsilon)} {2} 
   } = - \zeta(s)
 $$
   
  \noindent  That is  exactly  (\ref  {1p3}). %
  
  \noindent If we put $t \rightarrow t-i\epsilon$ in $i    \theta(t)=  i \Im \left[ \ln  \Gamma \left(    \frac{i t}{2} +\frac{1}{4}    \right )          -\frac{i t}{2}  \ln(\pi)  \right]  $ of  (\ref {ZMathem}), i.e. the exponential in  (\ref{Zt2}), we get:
  $$ i \left(   \beta\left[\frac{t-i\epsilon}{2} \ \ln\left( \frac{t-i\epsilon}{2 e}  \right) + \frac{\pi(\alpha +2)}{8} \right]
 -\frac{\beta(i t+\epsilon)}{2}\ln(\pi)  + O(t^{-1} ) \right) 
  = $$
  \be \label {intermediate1}
  i \left(   \beta \left[\frac{t-i\epsilon}{2} \ \left[  \ln \sqrt{\left( \frac{t^2 +\epsilon^2}{(2 e)^2}  \right) } + i \arctan\left( \frac{-\epsilon}{t}  \right)  \right] + \frac{\pi(\alpha +2)}{8} \right]
 -\frac{\beta(i t+\epsilon)}{2} \ln(\pi) + O(t^{-1} ) \right)
  \ee
  But:
   $$ \left[\frac{t-i\epsilon}{2} \ \left[  \ln \sqrt{\left( \frac{t^2 +\epsilon^2}{(2 e)^2}  \right) } + i \arctan\left( \frac{-\epsilon}{t}  \right)  \right] + \frac{\pi(\alpha +2)}{8} \right]
   =_{t \  big} - i 
    \frac{\epsilon}{2}
 \left( \ln\left ( \frac{t}{2e} \right)+ 1\right)+
 \frac{t}{2}
  \ln\left ( \frac{t}{2e} \right)+ \frac{\pi(\alpha +2)}{8}
   $$
\noindent  so (\ref{intermediate1}) becomes:

 $$=
   i \left(   \beta \left[   - i 
    \frac{\epsilon}{2}
 \left( \ln\left ( \frac{t}{2e} \right)+ 1\right)+
 \frac{t}{2}
  \ln\left ( \frac{t}{2e} \right)+ \frac{\pi(\alpha +2)}{8}      \right]
 -\frac{\beta(i t+\epsilon)}{2} \ln(\pi) + O(t^{-1} ) \right)=
  $$
   $$
      \left(   \beta \left[   +
    \frac{\epsilon}{2}
  \ln\left ( \frac{t}{2} \right)+i
 \frac{t}{2}
  \ln\left ( \frac{t}{2e} \right)+ i\frac{\pi(\alpha +2)}{8}      \right]
 -\frac{\beta(i t+\epsilon)}{2} \ln(\pi) + O(t^{-1} ) \right)
   $$

   \noindent 
    from   ( \ref {RealPart}) and  ( \ref  {ImmPart} ) we see  that, under $t \rightarrow t-i\epsilon$,   ( \ref {Zt2}) simplifies to:    %
   $$
   Z(t \rightarrow t- i \epsilon) =
   \frac {      e^{i\left[  \Im \left[ \ln  \Gamma \left(    \frac{i t +\epsilon}{2} +\frac{1}{4}    \right )      \right]     -\frac{t-i \epsilon}{2}  \ln(\pi)     \right]}                    (-\xi(s) )}  
   {e^{[ i \Im    \ln \Gamma(s/2) + \Re \ln  \Gamma(s/2)  ]} \left(  \frac{s(1-s)}{2}  \right) \pi^{\left(-\frac{1}{4}-\frac{\beta(i t)}{2} \right)}   } =
   $$
   $$
   =
     \frac  {   e^{ i \left(   \beta\left[\frac{t}{2} \ \ln\left( \frac{t}{2 e}  \right) + \frac{\pi(\alpha +2)}{8} \right]
       +
    \frac{\epsilon}{2} \ln\left ( \frac{t}{2} \right)
 -\frac{\beta(i t+\epsilon)}{2}\ln(\pi)  + O(t^{-1} ) \right)}  (- \xi(s))}  
   {e^{ \left [   +i \left[  \beta\left[\frac{t}{2} \ \ln\left( \frac{t}{2 e}  \right) + \frac{\pi(\alpha +2)}{8} \right]
 + \frac{\pi}{4} \epsilon  \right]+
   \left[  \frac{\alpha}{4} + \frac{1}{2} +   \frac{\beta \epsilon}{2}\right] \ln\left(\frac{t}{2} \right)
 +\ln(\sqrt{2 \pi  })-\frac{\pi}{4}t + O(t^{-1} )
    \right]}\left(  \frac{s(1-s)}{2} \right)\pi^{\left(-\frac{1}{4}-\frac{\beta(i t)}{2} \right)}  } =
   $$
   
   \be \label {TheSame}
 =   \frac {-\xi(s)}  
   {e^{[   +i  \frac{\pi}{4} \epsilon ]} \frac{s(1-s)}{2} 
\left( \frac{2}{ \pi t} \right)^{0.25} \sqrt{2 \pi  } 
e^{[   
 -\frac{\pi}{4}t + O(t^{-1} )
    ]} 
     }  
     \rightarrow 
      \frac {-\xi(s)}  
   {e^{[   +i  \frac{\pi}{4} \epsilon ]}       s(1-s)   
\left( \frac{ \pi}{2  t} \right)^{0.25}  
e^{[   
 -\frac{\pi}{4}t + O(t^{-1} )
    ]} 
     }  = Z(t,\epsilon)
   \ee

\noindent  So the almost identity of  numerical results  for    $Z(t -i\epsilon)$  and  $Z(t,\epsilon)$ in    \cite [p.~2]  {Giovanni Lodone 2021},
also  with   $ C_k(p,\epsilon) =0  \ ; \  \forall k  \ge1$ in (\ref  {RDue}), and $Err(t)=0$ 
in (\ref {ZSinhECosh}) ,

   \noindent is justified at $|t|$ big ( i.e. if  $t^{1/4} \approx [s(s-1)]^{1/8} $)    ) , where  (\ref  {TheSame}) holds.

\end{document}